%% file: mark.tex
\input amstex
\documentstyle {amsppt}
\tolerance=3000
\openup 6pt
\nologo

\topmatter
\title
Distortion results and invariant Cantor sets of unimodal maps\\
\endtitle

\author
Marco Martens
\endauthor

\affil
Institute for Mathematical Sciences\\
State University of New York at Stony Brook\\
Stony Brook USA
\endaffil


\abstract{A distortion theory is developed for $S-$unimodal maps. It will
be used to get some geometric understanding of invariant Cantor sets. In 
particular attracting Cantor sets turn out to have Lebesgue measure zero. 
Furthermore
the ergodic behavior of $S-$unimodal maps is classified according to a 
distortion property, called the Markov-property.}
\endabstract
\endtopmatter

\openup 6pt
\bigskip
\centerline{\bf 1. Introduction}
\bigskip

\flushpar
The work presented here originated in the question whether or not attracting
Cantor sets of unimodal maps have Lebesgue measure zero. This question led
to a general $S-$unimodal distortion theory. As applications of this theory 
we got uniform proofs of the basic known ergodic properties: ergodicity,
conservativity, existence of attractors. But also an answer to the original
question.

\proclaim{Theorem A} Cantor attractors of $S-$unimodal maps have Lebesgue 
measure zero.
\endproclaim

\flushpar
The strategy for studying invariant Cantor sets is constructing open,
arbitrarily fine  and nested covers of them. These covers are constructed in 
such a way that an invariance property appears: except for the component 
containing the critical point
every component is mapped monotonically onto its image which is also a 
component of the cover. Finally all components are transported to the central 
one, that is the one containing the critical point. The main question to be 
answered is whether this transport has good distortion properties. 
    
\flushpar
In fact the covers of the Cantor sets are part of covers of the almost the
whole interval,
having the same invariance property. $S-$unimodal maps having arbitrarily fine
covers with uniform good distortion properties are said to have the 
Weak-Markov-property.

\proclaim{Theorem B} Every $S-$unimodal map not having periodic attractors
has the Weak-Markov-Property.
\endproclaim

\flushpar
The tools developed for proving Theorem B are used to proof Theorem A. The 
ergodicity of $S-$unimodal maps is a direct consequence of the 
Weak-Markov-Property. For understanding stronger ergodic properties we need
a stronger distortion property.

\flushpar
Using the Weak-Markov-Property we see that smaller and smaller intervals are 
transported  with uniform bounded distortion to smaller and smaller central 
intervals. For getting stronger ergodic properties we need to find almost
everywhere
smaller and smaller intervals transported with uniform bounded distortion 
to a fixed big interval. Maps having this distortion property are said to 
have the Markov-Property.

\flushpar
Using this Markov-Property all basic ergodic properties can be proved in a
uniform way: conservativity, existence of attractors, etc. 
   
\flushpar
The main application of the  Markov-Property is an ergodic classification
of $S-$unimodal maps. In particular it can be used to classify the maps having
an attracting Cantor set in the sense of Milnor ([Mi]). 

\proclaim{Theorem C} An $S-$unimodal map not having a periodic attractor has 
the Markov-property if and only if it doesn't have a Cantor attractor.
\endproclaim

\flushpar
An appendix is added in which the basic notions of $S-$unimodal dynamics
are defined.

\flushpar
The results presented here are taken from the authors thesis defended at the
Technical University of Delft in 1990. 

\flushpar
The author would like to thank the Instituto de Matematica Pura e
Aplicada (IMPA) at Rio de Janeiro in which this work was done.
 
\input mark2.tex

\input mark3.tex
\input mark4.tex
\input mark5.tex

\input marka.tex
\input markr.tex

\bye

%% file: mark2.tex
\tolerance=3000

\bigskip
\centerline{\bf 2. Covers and induced maps}
\bigskip

\flushpar
The analytical properties of the construction presented here are based on the
following fundamental Lemma. Its proof can be found in different places, for 
example [MMS].

\proclaim{Lemma 2.1 (Koebe-Lemma)} Let $M,T$ be intervals in $[0,1]$ with $M\subset T$. The
components of $T\backslash M$ are denoted by $L$ and $R$. For every 
$\epsilon>0$ there exist $\delta>0$ and $K>0$ such that the following holds. 
Let $f:[0,1]\to [0,1]$ be a map with negative Schwarzian derivative. If $f^n|T$ 

is monotone and
$$
|f^n(L)|\ge \epsilon |f^n(M)| \text{ and } |f^n(R)|\ge \epsilon |f^n(M)| 
$$ 
then
\parindent=15pt
\item{1)} $\frac{|Df^n(x)|}
               {|Df^n(y)|}\le K$ for $x,y\in M$     \quad   (Koebe-Lemma);
\item{2)} $T$ contains a $\delta-$scaled neighborhood of $M$   
\quad    (Macroscopic-Koebe-Lemma).
\endproclaim

\proclaim{Corollary 2.2} For every $\rho>0$ there exists $\delta>0$ with the 
following property.  Let $\{M_i|i\ge 1\}$ be a pairwise disjoint collection
of subintervals of the interval $T$. For $i\ge 1$ let $L_i$ and $R_i$ be such 
that
\parindent=15pt
\item{1)} $L_i$ is a component of $T-\cup_{j\ge 1} M_j$ next to $M_i$; 
\item{2)} $R_i$ is a component of $T-M_i$ next to $M_i$;
\item{3)} $L_i\cap R_i=\emptyset$;
\item{4)} $|L_i|\ge \rho |M_i|$ and $R_i|\ge \rho |M_i|$.

\flushpar 
If $g:T'\to T$ is monotone, onto and has $Sg(x)<0$ for $x\in T'$ then for 
$i\ge 1$
$$
|L_i'|\ge \delta |M_i'| 
$$
where $L_i'$ and  $M_i'$ are the preimages under $g$ of respectively 
$L_i$ and  
$M_i$.  
\endproclaim

\flushpar
The Koebe-Lemma assures bounded distortion if there are big extensions of 
monotonicity on both sides. The next step is to develop topological 
instruments for studying  the maximal extensions.

\bigskip
\flushpar
In this section we will fix an $S-$unimodal map $f:[0,1]\to [0,1]$ with 
critical point $c$ and we assume that $f$ does not have periodic attractors.

\bigskip
\flushpar
For $x\in [0,1]$ denote the interval $(x,\tau(x))$ by $V_x$ (the involution 
$\tau$ is defined in the appendix). Furthermore define
$$
\Cal{N}=\{x\in [0,c)| V_x\cap orb(x)=\emptyset\}.
$$ 
The points in $\Cal N$ are called {\it nice}. Observe that every periodic orbit
contains nice points. Hence $\Cal N$ is not empty. Moreover since the critical 
point is accumulated by periodic orbits $\Cal N$ also accumulates on $c$. 
Clearly $\Cal N$ is closed. Fix $x\in \Cal N$.

\proclaim{Lemma 2.3} For $i=1,2$ let $T_i\subset [0,1]$ be two different
intervals such that 
$f^{n_i}:T_i\to V_x$ is monotone and onto for some $n_1\le n_2$ . If 
$T_1\cap T_2\ne \emptyset$ then $T_2\subset T_1$ and $n_1<n_2$.
\endproclaim

\demo{proof} Suppose that $\partial T_1\cap T_2\ne \emptyset$. Then 
$n_1\ne n_2$. And $x\in f^{n_1}(T_2)$ which implies $f^{n_2-n_1}(x)\in V_x$. Contradiction.
\qed
\enddemo

\flushpar
The set of points who visit $V_x$ is called $C_x\subset [0,1]$. Observe that 
it can be described as being the union of all intervals $T$ such that 
$f^n:T\to V_x$ is monotone and onto for some $n\ge 0$. Furthermore let 
$\Lambda_x=[0,1]-C_x$ and $D_x=f^{-1}(C_x)\cap V_x$. 

\proclaim{Lemma 2.4} Let $I\subset C_x$ be a component. Then there exists 
$n\ge 0$ such that $f^n:I\to V_x$ is monotone and onto. 
Furthermore $\{I,f(I),\dots,f^n(I)=V_x\}$ is a pairwise disjoint collection.
In particular there exists only one such $n\ge 0$.
\endproclaim

\demo{proof}
Lemma 2.3 implies easily that for every component $I$ of $C_x$ there exists 
$n\ge 0$ such that $f^n:I\to V_x$ is monotone and onto. Observe that this 
number $n\ge 0$ is defined uniquely: $c\in f^n(I)$ which implies that 
$f^{n+s}|I$ is not monotone for $s>0$. 

\flushpar
To proof the disjointness of the orbit it suffices to proof that  
$f^j(I)\cap V_x=\emptyset$ for all $j<n$. Suppose $f^j(I)\cap V_x\ne \emptyset$
with $j<n$. By Lemma 2.3 we get $f^j(I)\subset V_x$.

\flushpar
Let $j<n$ be minimal such that $f^j(I)\subset V_x$. Let $H\subset [0,1]$ be the
maximal interval containing $I$ with $f^j|H$ is monotone and 
$f^j(H)\subset V_x$. Suppose that a component $L$ of $H-I$
is mapped in $V_x$: $\overline{f^n(L)}\subset V_x$. Then by maximality there 
exists $s<j$ such that $c\in \partial f^s(L)$. The minimality of $j$ 
implies $f^s(I)\cap V_x=\emptyset$. Hence $x\in f^s(L)$. Because 
$\overline{f^j(L)}\subset V_x$ this implies $orb(x)\cap V_x\ne \emptyset$.
Contradiction.

\flushpar
So we conclude that $f^j:H\to V_x$ is monotone and onto. Furthermore $H$
contains the component $I$ of $C_x$. The definition of $C_x$ implies that
$H=I$. In particular $j=n$. Contradiction. 
\qed
\enddemo

\flushpar
Lemma 2.4 states the invariance property of the covers $C_x$ discussed in
the introduction. It enables us to define the {\it Transfer map}
$$
T_x:C_x\to V_x
$$ 
and also the {\it Poincar\'e map}
$$
R_x:D_x\to V_x
$$ 
by $R_x=T_x\circ f$.

\flushpar
The next lemma shows that these induced maps are defined almost everywhere.

\proclaim{Lemma 2.5} If $x\in \Cal{N}$ then 
\parindent=15pt
\item{1)} $|C_x|=1$ (and $|\Lambda_x|=0)$;
\item{2)} $|D_x|=|V_x|$.

\flushpar
Furthermore 
\parindent=15pt
\item{3)} $\Lambda_x$ is invariant;
\item{4)} if $y\in \Lambda_x$ is such that $orb(y)\cap \overline V_x=\emptyset$
then $\Lambda_x$ accumulates from both sides on $y$.
\endproclaim

\demo{proof}
Take $y\in \Lambda_x$ and assume $f(y)\notin \Lambda_x$, say 
$f(y)\in I\in C_x$. If $c_1\notin I$ then $f^{-1}(I)\subset C_x$. Hence 
$c_1\in I$. Because $x\in \Cal N$ we get easily that $f^{-1}(I)\subset V_x$.
This gives the contradiction $y\in V_x$. So $\Lambda_x$ is invariant.

\flushpar
To prove the first two statements it suffices to proof $|\Lambda_x|=0$. 
Because $\Lambda_x$ is a closed 0-dimensional invariant set not 
containing the critical point a 
well-known lemma in [Mi] states that $\Lambda_x$ has Lebesgue measure zero.

\flushpar
Take $y\in \Lambda_x$ with $orb(y)\cap \overline{V_x}=\emptyset$. Suppose that
$y$ is not accumulated from both sides by $\Lambda_x$. Hence $y$ is a boundary 
point of a component of $C_x$. Lemma 2.4 states that the orbit of $y$ passes 
through the boundary of $V_x$. Contradiction.
\qed
\enddemo

\flushpar
As a consequence of the previous lemma the set
$$
\Cal{R}=\{x\in [0,1]|c\in \omega(x)\}
$$
has Lebesgue measure $1$.

\bigskip
\flushpar
The branches of the Transfer map are monotone. To apply the Koebe-Lemma we have
to know how much the monotonicity can be extended.

\flushpar
Suppose $c_1\in C_x$. So there exists a component $S_x\subset C_x$ with
$c_1\in S_x$. Let $\psi(x)=\partial f^{-1}(S_x)\cap [0,c)$ and define 
$U_x=V_{\psi(x)}=f^{-1}(S_x)$. 

\flushpar
We get $S_x\subset (f(V_x),1)$ because $x$ is nice. Hence 
$U_x=V_{\psi(x)}\subset V_x$. Using lemma 2.4 we get 
$orb(f(\psi(x)))\cap V_x=\emptyset$. In particular $\psi(x)\in \Cal N$.

\flushpar
To finish the definition of the increasing function $\psi:\Cal N\to \Cal N$
let $\psi(x)=x$ if $c_1\ne C_x$.

\flushpar
The pair $(V_x,U_x)$ is called a {\it transfer range}. If $V_x$ contains a
$\delta-$scaled neighborhood of $U_x$ then the pair is called a 
$\delta-$transfer range.

\proclaim{Proposition 2.6} Let $x\in \Cal N$ and $I\subset C_{\psi(x)}$ be a 
component, say $T_{\psi(x)}|I=f^n|I$. Then there exists an interval $T_I$
containing $I$ such that 
$$
f^n:T_I\to V_x
$$
is monotone and onto.
\endproclaim

\flushpar
Combining the Koebe-Lemma with this proposition we get that the branches of 
$T_{\psi(x)}:C_{\psi(x)}\to U_x$ have uniformly bounded distortion. The bound
just depends on  the space in $V_x$ around $U_x$.

\flushpar
A pair $(T,I)$ of intervals with the property that for some $n\ge 0$, called
the {\it transfer time},
\parindent=15pt
\item{1)} $f^n:T\to V_x$ is monotone and onto;
\item{2)} $f^n(I)=U_x$

\flushpar
is called {\it TI-pair} for $(V_x,U_x)$. As in lemma 2.4 observe that the time
$n\ge 0$ is defined uniquely. As we saw above every TI-pair for a certain 
transfer range has uniform bounded distortion on the middle part. 

\demo{proof of proposition 2.6} Let $I$ be a component of $C_{\psi(x)}$, say 
$T_{\psi(x)}|I=f^n|I$, and let $H\subset [0,1]$ be the maximal interval 
containing $I$ with $f^n|H$ is monotone and $f^n(H)\subset V_x$. Suppose by 
contradiction that the component $L$ of $H-I$ is mapped in 
$V_x$: $\overline{f^n(L)}\subset V_x$. Then by maximality there exists $j<n$ 
such that $c\in \partial f^j(L)$. We know from lemma 2.4 that 
$f^j(I)\cap U_x=\emptyset$. Hence $\psi(x)\in f^j(L)$. Because 
$\overline{f^n(L)}\subset V_x$ this implies $orb(\psi(x))\cap V_x\ne \emptyset$.
Contradiction.
\qed
\enddemo

\flushpar
The intersection behavior of the TI-pairs is formulated in

\proclaim{Lemma 2.7} For $i=1,2$ let $(T_i,I_i)$ be TI-pairs for $(V_x,U_x)$ with
transfer times respectively $n_1$ and $n_2$. If $T_2\cap T_1\ne \emptyset$ and 
$n_2\ge n_1$ then
$$
T_2\subset T_1-I_1  \text{ or } T_2\subset I_1
$$ 
and $n_2>n_1$.
\endproclaim    

\demo{proof}
From lemma 2.3 we get $T_2\subset T_1$ and $n_2<n_1$. Suppose, by 
contradiction, $T_2\cap \partial I_1\ne \emptyset$. Then $\psi(x)\in f^{n_1}(T_2)$. 
Hence $f^{n_2-n_1}(\psi(x))\in V_x$. Contradiction.  
\qed
\enddemo

\flushpar
We finish with the definition of the last type of induced maps. Let $E_x$
be the union of all intervals $I$ which are part of a TI-pair $(T_I,I)$
for $(V_x,U_x)$ with {\it positive} transfer time.  

\flushpar
Again lemma 2.7 implies that the connected components of $E_x$ are intervals 
which are part of a TI-pair for $(V_x,U_x)$. This defines the {\it Markov} map
$$
M_x:E_x\to U_x.
$$
Directly after defining the first two induced maps we where able to show that 
they are defined almost everywhere. For Markov maps the situation is not that
simple. Section 4 will deal with the question whether or not the Markov maps
are defined almost everywhere. There we will show that Markov-maps are 
defined almost everywhere only in the case of absence of Cantor attractors.

%% file: mark3.tex
\tolerance=3000

\bigskip
\centerline{\bf 3. The Weak-Markov-Property}
\bigskip

\flushpar
In this section we are going to prove a general distortion result for 
$S-$unimodal maps, called the weak-Markov-property. 

\proclaim{Definition 3.1} An $S-$unimodal map is said to satisfy the 
{\it Weak-Markov-Property} if there exist $K>0$, a set 
$\Cal{G}^{w}\subset [0,1]$ with full Lebesgue measure and a set 
$\Cal D\subset \Cal N$ accumulating at the critical point such that for every
$x\in \Cal{G}^{w}$ and for every $y\in \Cal D$ there exist $t\ge 0$ and an 
interval $I\ni x$ such that 
$$
f^t:I\to V_y
$$
is a diffeomorphism with distortion bounded by $K$.
\endproclaim

\proclaim{Theorem 3.2} Let $f$ be an unimodal map
without periodic attractor. There exists $\delta>0$ and $\Cal D \subset \Cal N$
accumulating at the critical point such that for every $x\in \Cal D$ the 
following property holds.

\flushpar
If $I$ is a component of $C_x$, say $T_x|I=f^n|I$ then $f^n$ maps the maximal 
interval containing $I$ on which $f^n$ is monotone over a $\delta-$scaled 
neighborhood of $V_x$.  
\endproclaim

\flushpar
The Koebe-lemma implies directly the main consequence of this Theorem: 

\proclaim{Theorem 3.3 (The Weak-Markov-Property)} 
Every $S-$unimodal map without periodic attractor satisfies the 
Weak-Markov-Property. 
\endproclaim

\flushpar
During the proof we will see that the number $\delta>0$ only depends on the 
power of the critical point 
and that the set $\Cal D$ has a topological definition. 

\flushpar
During the proof of Theorem 3.3 we will see that in general the monotone 
extensions of the branches $f^n:I\to U_x$ with $I\subset C_x$ will have images
much bigger than the transfer range. So Theorem 3.3 does not learn us 
something about the geometry of transfer ranges.
 
\flushpar
However in our study of the Lebesgue measure of 
the critical limit sets we need a better understanding of the geometry of 
transfer ranges. That is why we need the following stronger form of 
Theorem 3.3.

\proclaim{Theorem 3.4} Let $f$ be an only finitely renormalizable unimodal map 
having a non-periodic recurrent critical point. Then there exists a sequence  
$\{(U_n,V_n)\}_{n\ge 0}$  of $\delta-$transfer ranges with $|V_n|\to 0$. 
\endproclaim

\flushpar
In this section we are going to prove Theorem 3.2 and 3.4. Fix an $S-$unimodal
map $f$ whose  critical point $c$ is
recurrent. Furthermore we assume $f$ not to have periodic attractors.

\proclaim{Proposition 3.5} There exist $\delta,\rho>0$ such that for 
$x\in \Cal{N}$ the following holds.
\parindent=15pt
\item{1)} Assume $c\in R_x(U_x)$. Let $I\subset C_{\psi(x)}$ be a component, 
say 
$T_{\psi(x)}|I=f^k|I$, and $T$ the maximal interval containing $I$ for which 
$f^k$ is monotone. Then $f^k(T)$ contains a $\delta-$scaled neighborhood of 
$U_x$. 
\item{2)} Assume $c\notin R_x(U_x)$ and $|V_x|\le (1+\rho)|U_x|$. If 
$T_x|S_x=f^n|S_x$ and if $T$ is the maximal interval containing $S_x$ such
that $f^n|T$ is monotone then $f^n(T)$ 
contains a $\delta-$scaled neighborhood of $[R_x(U_x),c]$.  
\endproclaim

\demo{proof} Let $T_x|S_x=f^n|S_x$ and $M=f(U_x)$. We are going to study the orbit
$\{M,f(M),...,f^n(M)\}$. $M$ is contained in $S_x$. Hence by lemma 2.4 this 
orbit consists
of pairwise disjoint intervals. Choose $m\in \{0,1,2,...,n\}$  such that 
$|f^m(M)|\le |Df|_0^2|f^j(M)|$ for all $j\le n$ and $f^m(M)$ has neighbors 
on both sides. This means that both components of $[0,1]-f^m(M)$ contains 
intervals of the form $f^i(M)$ with $i\le n$. Observe that by taking $m'$ in 
such a way that $f^{m'}(M)$ is the smallest one we can take 
$m\in \{m',m'+1,m'+2\}$  having the described property. Let $f^l(M)$ and 
$f^r(M)$ be the direct neighbors of $f^m(M)$. 

\flushpar
Let $H$ be the maximal interval containing $M$ for which $f^m|H$ is monotone
and $f^m(H)\subset [f^l(M),f^r(M)]$. We claim
$$
f^m(H)=[f^l(M),f^r(M)].
$$ 
Fix a component $L$ of $H-M$ and assume that $\overline{f^m(L)}\subset
[f^l(M),f^r(M)]$. by maximality there exists $j\le n$ with 
$c\in \partial{f^j(L)}$. Because $f^j(M)\cap U_x=\emptyset$ we see that 
$f^{m-j-1}(M)\subset \overline{f^m(L)}\subset [f^l(M),f^r(M)]$. But $m-j-1\le
m-1$. Hence $[f^l(M),f^r(M)]$ contains at least four intervals of the form
$f^i(M)$, $i\le n$. This contradiction finish the proof of the claim. 

\flushpar
From the definition of $m,l,r$ and the Macroscopic-Koebe-lemma we get a 
universal constant $\delta_1>0$ (only depending on $|Df|_0$) such that $H$ 
contains a $\delta_1-$scaled neighborhood of $M$. Because the critical point 
of $f$ is non-flat there exists a universal constant $\delta_2>0$ such that 
$$
H'=f^{-1}(H) \text{ contains a }
             \delta_2-
             \text{scaled neighborhood of }
             U_x.
$$

\demo{proof of statement 1} Assume that $c\in f^n(M)$. Take a component $I$ of
$D_{\psi(x)}$, say with transfer time  $k\ge 0$. Denote by $T$ the maximal interval 
containing $I$ on which $f^k$ is monotone and $f^k(T)\subset H'$. We claim
$f^k(T)=H'$ which proves statement 1.

\flushpar
To prove this, fix a component $L$ of $T-I$ and suppose that $\overline{f^k(L)}
\subset H'$. From the maximality of $T$ we get $j< k$  with $c\in \partial 
f^j(L)$. But $f^j(I)\cap U_x=\emptyset$. Therefore $f^k(L)$ contains 
$f^{k-j-1}(M)$ in its closure which is contained in $H'$. Because $f^{m+1}$
is monotone on the component of $H'-\{c\}$ which contains $f^k(L)$, the iterate 
$f^{m+1+k-j-1}$ maps the closure of $M$ monotonically
into  the interior of $[f^l(M),f^r(M)]$. In particular, because $c\in f^n(M)$,
we have  $k+m-j\le n$. Again this implies that  $[f^l(M),f^r(M)]$ contains
at least four intervals of the form $f^i(M), i\le n$. Contradiction.
\enddemo

\demo{proof of statement 2} Let $T$ be the maximal interval containing $S_x$ 
on which $f^n$ is  monotone. The components of $T-M$ are denoted by $L,R$, say
$c\in f^n(R)$. Because $f^n(S_x)=V_x$ and $c\notin f^n(M)$ the interval 
$f^n(R)$ contains a component of $V_x-\{c\}$. So using the non-flatness of the
critical point we get a constant $\delta_3>0$ such that
$$
|f^n(R)|\ge \delta_3 |[f^n(M),c]|.
$$  
So we only have to study $f^n(L)$. We claim
$$
L'\subset f^n(L)
$$
where $L'$ is the component of $H'-U_x$ which lies on the same side of $c$ as 
$f^n(L)$. Once this claim is proved statement 2 follows by taking $\rho>0$ 
sufficiently small.

\flushpar
Indeed, suppose by contradiction $\overline{f^n(L)} \subset L'$. Again the 
maximality of $T$ assures the existence of $j< n$ with $c\in \partial
f^j(L)$. Observe $f^j(M) \cap V_x=\emptyset$. So we get $f^{n-j-1}(M)\subset 
\overline{f^n(L)} \subset L'$. Consider the interval $[f^j(M),c]$ (if 
$f^j(M)$ and $f^{n-j-1}(M)$ are on the same side of $c$. Otherwise consider 
$\tau([f^j(M),c])$.

\flushpar
The map $f^{n-j}: [f^j(M),c] \to [f^n(M),f^{n-j-1}(M)]$ is monotone and onto. 
Because $f$ doesn't have a periodic attractor we get
$$
f^j(M) \subset (f^{n-j-1}(M),f^n(M)).
$$
Because $f^{n-j-1}(M)\subset L'$ and hence $f^{n-j-1+m+1}(M) \subset 
[f^l(M),f^r(M)]$ we get $n-j-1+m+1>n$. Hence
$$
j<m.
$$
Furthermore $[f^{n-j-1}(M),c] \subset [L',c]$ and $f^{m+1}|[L',c]$ is 
monotone. Because $j<m$ the map $f^{j+1}:[f^{n-j-1}(M),c]\to [f^j(M), f^n(M)]$ 
is monotone and onto. Furthermore we have $[f^j(M),f^n(M)]\subset 
[f^{n-j-1}(M),c]$. Hence $f$ has a periodic attractor. Contradiction.
\qed
\enddemo
\enddemo

\flushpar
Shortly speaking Proposition 3.5 states that the central branch of the 
Poincar\'e map $R_x$ has a quadratic shape. 

\proclaim{Corollary 3.6} Let $f$ be an unimodal not having periodic attractors.
There exist $\rho>0$ and $K<\infty$ with the following properties. Let 
$x\in \Cal N$ and suppose $T_x|S_x=f^n$. If $|V_x|\le (1+\rho)|U_x|$ then
\parindent=15pt
\item{1)} For all $y_1,y_2\in f(U_x)$
$$
\frac1K\le \frac{|Df^n(y_1)|}{|Df^n(y_2)|} \le K;
$$
\item{2)} For all $y\in U_x$
$$
|Df^{n+1}(y)|\le K.
$$
\endproclaim

\demo{proof} The first statement follows directly by applying the Koebe-lemma
and the previous Proposition.

\flushpar
To prove the second statement we define a point $m\in \overline{U_x}$. Let $m=\psi(x)$ if 
$|Df^{n+1}(\psi(x))|\le \frac{3}{1-\rho}$. If $|Df^{n+1}(\psi(x))|> 
\frac{3}{1-\rho}$ let $m$ be the closest point in $U_x$ to $\psi(x)$ such that
$|Df^{n+1}(m)|=\frac{3}{1-\rho}$ (because $Df^{n+1}(c)=0$ such an $m$ exists).

\flushpar
First we are going to show 
$$
\frac{|[m,c]|}{|[\psi(x),c]|}\ge \frac13.
$$ 
We may assume $m\ne \psi(x)$. Because $f^{n+1}(U_x) \subset V_x$ we have
$$
|V_x|\ge |f^{n+1}([\psi(x),m])|\ge \frac{3}{1-\rho}|[\psi(x),m]|.
$$
Hence
$$
\align
\frac{|[m,c]|}{|[\psi(x),c]|}
&=   \frac{|[\psi(x),c]|-|[\psi(x),m]|}{|[\psi(x),c]|}\ge\\
&\ge 1-2\frac{1-\rho}{3}\frac{|V_x|}{|U_x|}\ge 1-2\frac{1-\rho}{3}
                                                  \frac1{1-\rho}=\frac13. 
\endalign
$$
Because $c$ is non-flat we get for all $y\in U_x$
$$
\align
|Df^{n+1}(y)|&=\frac{|Df^{n+1}(y)|}{|Df^{n+1}(m)|}. |Df^{n+1}(m)|\le
               \frac{|Df(y)|}{|Df(m)|}K\frac{3}{1-\rho}\le\\
             &\le B\frac{|x-c|^\alpha}{|m-c|^\alpha}\le 
                  B\frac{|\psi(x)-c|^\alpha}{|m-c|^\alpha}
              \le B 3^\alpha=A
\endalign
$$
where $\alpha$ is the order of the critical point and $B$ a positive constant
depending on $\rho, K$ and the behavior of $f$ around the critical point.
\qed
\enddemo

\proclaim{Lemma 3.7} There exists $\delta>0$  such that for all $x\in \Cal N$
with
\parindent=15pt
\item{1)} $c\notin R_x(U_x)$;
\item{2)} $R_x(c)\in V_x-U_x$

\flushpar
$V_x$ contains a $\delta-$scaled neighborhood of $U_x$.
\endproclaim

\demo{proof}
Let $\delta>0$ be given by proposition 3.5 and denote the part of the 
$\frac{\delta}2-$scaled neighborhood of $[R_x(U_x),c]$ which lies on the same
side of $c$ as $R_x(U_x)$ by $H'$. The pair of intervals $P\subset Q$ is such 
that
\parindent=15pt
\item{1)} $f(U_x)\subset P$ and $P\subset S_x$;
\item{2)} $f^n(P)=[R_x(U_x),c]$;
\item{3)} $f^n:Q\to [H',c)$ is monotone and  onto (where 
$R_x|U_x=f^{n+1}|U_x$);

\flushpar
Proposition 3.5 implies that such intervals exist. Let $|V_x|=(1+\rho)|U_x|$
and assume that $\rho$ is small. Furthermore the Koebe-Lemma
tells us that $f^n|Q$ has bounded distortion.  We 
have to show that $\rho$ is not too small. Observe that $|R_x(U_x)|=O(\rho)$.
This follows from assumption 1) and 2). 

\flushpar
Consider $U=f^{-1}(Q)$ as a neighborhood of $U_x$. Using the bounded 
distortion of $f^n|Q$ and the fact that the critical point of $f$ is of order
$\alpha$ it is easy to show that $U$ is a $O(\rho^{-\frac1\alpha})-$scaled 
neighborhood of $U_x$. This bound implies that for $\rho$ small $H'\subset U$.

\flushpar
Suppose that $\rho$ is sufficiently small to assure that $U$ will contain 
$H'$. Now let $H$ be the component of $U-\{c\}$ containing $R_x(U_x)$. Then 
$$
f^{n+1}:H\to H'\subset H
$$
is monotone. Hence $f$ has a periodic attractor. Contradiction: $\rho$ is away
from zero. 
\qed
\enddemo

\proclaim{Lemma 3.8} There exist $\rho,\delta>0$ such that for all 
$x\in \Cal N$ with
\parindent=15pt
\item{1)} $|V_x|\le (1+\rho)|U_x|$;
\item{2)} there exists $p\in U_x$ with $R_x(p)=p$ and $V_p\subset R_x(V_p)$

\flushpar
$V_p$ contains a $\delta-$scaled neighborhood of $U_p$.
\endproclaim

\demo{proof}
Observe that this periodic point is nice: $p\in \Cal{N}$. Let $M\subset V_p$ be
the maximal interval with $R_x(M)\cap V_p=\emptyset$. Easily we get $U_p\subset M$. 
Let $\rho$ be given by corollary 3.6. Then we get a universal bound on 
$DR_x|U_x$: $|DR_x(y)|\le A$ for all $y\in V_p$. This implies that both 
components of $V_p-M$ are bigger than $\frac1A|V_p|$. Hence
$$
|U_p|\le |M| \le (1-\frac2A)|V_p|.
$$ 
Let $\delta=\frac1A$.
\qed
\enddemo

\demo{proof of Theorem 3.2}
If there exists a neighborhood $V$ of $c$ such that 
$\overline{V}\cap \overline{orb\{c\}}=\emptyset$ then let $\Cal D = \Cal N \cap V$. The 
theorem follows from the fact that the boundary of the images of the maximal 
intervals of monotonicity consists of critical values, points in the orbit
of $c$. Hence we may assume that the critical point is recurrent.

\flushpar
First we are going to define the sequence of closest approach to $c$. Let 
$c_n=f^n(c)$ and define
$$
\align
q(1)&=1;\\
q(n+1)&=\min\{t\in \bold{N}| c_t\in V_{c_{q(n)}}\}.
\endalign
$$
Because $c\in \omega(c)$ and $c$ is not periodic the sequence 
$\{q(n)\}_{n\ge 0}$ is well defined. Since $c$ is an accumulation point of 
$\Cal N$ there are infinitely many $n\ge 1$ for which 
$(V_{c_{q(n-1)}}-V_{c_{q(n)}})\cap \Cal N \ne \emptyset$. For those $n\ge 1$ 
define
$$
x(n)=\sup\{(V_{c_{q(n-1)}}-V_{c_{q(n)}})\cap \Cal N \cap [0,c)\}.
$$
Observe that the points $c_{q(n)}$ are not in $\Cal N$: the supremum is in fact
a maximum. Furthermore $x(n)\in (V_{c_q(n-1)}-V_{c_{q(n)}}\cap \Cal N)$ .
Hence for $i<q(n)$ 
$f^i(c)\notin V_{x(n)}$ and $f^{q(n)}(c) \in V_{c_{q(n)}} \subset V_{x(n)}$. In 
particular $R_{x(n)}|U_{x(n)}=f^{q(n)}$. 

\flushpar
Let $\Cal D = \{\psi(x(n))\}$. We distinguish two cases.

\demo{Low Case. $c\notin R_{x(n)}(U_{x(n)})$} 
Because $f$ doesn't have periodic attractors and $c\notin R_{x(n)}(U_{x(n)})$ 
we 
get $\psi(x(n))>x(n)$. Hence, by using the definition of $x(n)$ we get 
$\psi(x(n))\in V_{c_{q(n)}}$ or $U_{x(n)}\subset V_{q(n)}$.
So we can apply Lemma 3.7: $V_{x(n)}$ is a $\delta-$scaled neighborhood of
$U_{x(n)}$. Use proposition 2.6, the fundamental property of transfer ranges,
to finish the proof.
\enddemo

\demo{High Case. $c\in R_{x(n)}(U_{x(n)})$}
The statement follows directly from proposition 3.5.
\qed
\enddemo
\enddemo

\bigskip
\demo{proof of theorem 3.4}
The transfer ranges will be defined in terms of the points $y(n)\in \Cal N$.  
These points are an adjustment of the points $x(n)$ from the proof of theorem 
3.2 and will be defined below. Let $\delta>0$ be the minimum of the numbers 
$\delta$ given by lemma 3.7 and 3.8.  We distinguish two cases. 

\demo{Low Case. $c\notin R_{x(n)}(U_{x(n)})$} 
The definition of $x(n)$ implies that $R_{x(n)}(c)=c_{q(n)}\notin U_{x(n)}$. 
Hence we can apply Lemma 3.7: $V_{x(n)}$ is a $\delta-$scaled neighborhood of
$U_{x(n)}$. Let $y(n)=x(n)$. 
\enddemo

\demo{High Case. $c\in R_{x(n)}(U_{x(n)})$}
Let $\rho>0$ be given by lemma 3.8. If $|V_{x(n)}|\ge (1+\rho)|U_{x(n)}|$ then
let $y(n)=x(n)$. Assume that this metrical property doesn't hold. The map 
$R_{x(n)}|U_{x(n)}$ has a periodic point $p\in U_{x(n)}$. Because $f$ is only 
finitely renormalizable we may assume (by taking $n$ large enough) that 
$V_p\subset R_{x(n)}(V_p)$. Let $y(n)=p$ and apply lemma 3.8.   
\enddemo

\flushpar
Because $c_{q(n)}\to c$ it follows that $y(n)\to c$. 
\qed
\enddemo   

\flushpar
Let us finish this section with a simple consequence of the 
Weak-Markov-Property. It implies directly the ergodicity.

\proclaim{Lemma 3.9} Let $f$ be a $S-$unimodal map without periodic attractor
and $\Cal D\subset [0,1]$ the set given by the Weak-Markov-Property
(Theorem 3.3).
 If $X\subset [0,1]$ is
an invariant set with positive Lebesgue measure then
$$
\lim_{\Cal D \ni x\to c} \frac{|X\cap V_x|}{|V_x|}=1.
$$ 
\endproclaim 

\demo{Proof} Using lemma 2.5(1) and $|X|>0$ we get a density point $y\in X$ of $X$
with $y\in C_x$ for $x\in \Cal D$. For every $x\in \Cal D$ let $I_x$ be the 
component of $C_x$ containing $y$. Say $f^{n_x}:I_x\to V_x$ is diffeomorphic with 
distortion bounded by $K$. Furthermore the contraction principle implies
$|I_x|\to 0$ if $\Cal D \ni x\to c$. 
Consider
$$
\aligned
\lim_{\Cal{D}\ni x\to c} \frac{|X^c\cap V_x|}{|V_x|}&\le
\lim_{\Cal{D}\ni x\to c} \frac{|T_x(X^c\cap I_x)|}{|T_x(I_x)|}\\
&=\lim_{\Cal{D}\ni x\to c}\frac{\int_{X^c\cap I_x} |DT_x(t)|dt}
                              {|DT_x(\beta_n)||I_x|}  
\endaligned
$$
for some $\beta_n\in I_n$. Thus
$$
\lim_{\Cal D \ni x \to c} \frac{|X^c\cap V_x|}{|V_x|}
\le \lim_{\Cal D \ni x \to c} K \frac{|X^c\cap I_x|}{|I_x|}=0.
$$     
This finishes the proof. 
\qed
\enddemo

%% file: mark4.tex
\tolerance=3000

\bigskip
\centerline{\bf 4. The Markov-Property}
\bigskip

\flushpar
In this section we will study a second distortion result, called the
{\it Markov-Property}. This property is much stronger than the 
Weak-Markov-Property and serves for studying more complicated ergodic
properties. 

\proclaim{Definition 4.1} An $S-$unimodal map is said to satisfy the 
{\it Markov-Property} if there exists a set $\Cal{G}\subset [0,1]$ with full 
Lebesgue measure and $y\in \Cal N$ such that 
for every $x\in \Cal{G}$ there exist a sequence of TI-pairs $(T_n,I_n)$ for 
the transfer range $(V_y,U_y)$ with $x\in I_n$ and $|I_n|\to 0$.

\flushpar
In particular, if $U_y\ne V_y$, the maps 
$$
I_n\to U_y
$$ 
are diffeomorphisms with uniform  bounded distortion.
\endproclaim
 
\flushpar
We can reformulate the Markov-Property in terms of Markov-maps:
   
\proclaim{Lemma 4.2} An unimodal map $f$ has the Markov-Property
iff there exists $x\in \Cal{N}$ such that the Markov map $M_x:E_x\to U_x$ has 
full domain:
$$
|E_x|=1.
$$
\endproclaim

\proclaim{Proposition 4.3} Let $f$ be an $S-$unimodal map without periodic 
attractor. If the limit set of the critical point is not minimal 
then $f$ has the Markov-property.
\endproclaim   

\flushpar
Let us prove this proposition for the $S-$unimodal map $f$ satisfying the two 
conditions of the proposition. In particular $f$ is not infinitely 
renormalisable. 

\flushpar
Let us first deal with the case when $f$ is Misiurewicz, $c\notin \omega(c)$.
Take $x\in \Cal N$ such that $\overline{orb(c_1)}\cap V_x=\emptyset$. Then 
$c_1\notin C_x$ and $\psi(x)=x$. In this situation it is easy to show that
$$
E_x=(C_x-V_x)\cup (f^{-1}(C_x)\cap V_x).
$$ 
This implies, by using lemma 2.5, that $M_x$ has full domain: $|E_x|=1$.

\flushpar
In the sequel we will assume that the critical point $c$ is recurrent.
 
\proclaim{Lemma 4.4} There exists $x\in \Cal N$ and a sequence of intervals 
$K_n$, $n\ge 1$ such that for $n\ge 1$
\parindent=15pt
\item{1)} There are no renormalisations possible in $V_x$. In particular
$\psi(x)>x$;
\item{2)} $\partial K_n\subset \Lambda_x$ and $K_n\cap V_x=\emptyset$;
\item{3)} $K_n\cap \omega(c)\ne \emptyset$;
\item{4)} $|K_n|\to 0$.
\endproclaim

\demo{proof}
Let $q\in\omega(c)$ be such that $c\notin\omega(q)$. This is possible because
$\omega(c)$ is not minimal. Choose $x\in 
\Cal N$ such that
$$
orb(q)\cap \overline{V_x}=\emptyset.
$$ 
and such that there are no renormalisations possible in $V_x$.
Then $q\in \Lambda_x$ and $q$ is accumulated from both sides by $\Lambda_x$
(see lemma 2.5). Hence we can take a sequence of intervals $q\in K_n$ with
$|K_n|\to 0$ and $\partial K_n\subset \Lambda_x$ accumulating at 
$q\in\omega(c)$.
\qed
\enddemo

\flushpar
Let $x\in \Cal N$ and the intervals $K_n$ be given by lemma 4.4.
Because $K_n \cap \omega(c)\ne \emptyset$ there exists a sequence 
$t_n\to \infty$ with $c_j\notin K_n$ for $j<t_n$ and $c_{t_n}\in K_n$. We 
claim that for every $n\ge 1$ there exists an interval $K'_n$ containing $c_1$
such that
$$
f^{t_n-1}:K'_n\to K_n
$$ 
is monotone and onto. Indeed let $K'_n$ be the maximal interval containing 
$c_1$ with $f^{t_n-1}|K'_n$ is monotone and $f^{t_n-1}(K'_n)\subset K_n$. 
Assume by contradiction that $\overline{f^{t_n-1}(L'_n)}\subset K_n$ where 
$L'_n$ is a component of $K'_n-\{c_1\}$. Then the
maximality of $K'_n$ implies the existence of $j<t_n-1$ such that 
$c\in \partial f^j(L'_n)$. So $f^{t_n-1-j}(c)\in K_n$. Contradicting the fact
that $c_{t_n}$ is the first critical value in $K_n$. This proves the claim.

\flushpar
We are going to use the components of $C_{\psi(x)}$, which belongs to TI-pairs,
in $K_n$ to show that the domain of the Markov map $M_x$ is has positive upper 
density in $c$. This is done by pulling back  the TI-pairs into $K'_n$, close 
to $c_1$. And then one step more to get them close to $c$. However we have 
to be careful. The whole TI-pairs have to be pulled back. The collection 
$\hat{\Cal{M}}_n$ defined below will contain the components for which this is 
impossible.

\flushpar
The statements in the lemma below follow directly from respectively
 proposition 2.6,
the definition of $C_x$ and lemma 2.5. We will leave the the proofs
to the reader.

\proclaim{Lemma 4.5} If $I$ is a component of $C_{\psi(x)}$ with 
$I\cap K_n\ne \emptyset$ then there exists an interval $T$ containing $I$
such that 
\parindent=15pt
\item{1)} $(T,I)$ is a TI-pair for $(V_x,U_x)$;
\item{2)} $T\subset K_n$.

\flushpar
Furthermore 
$$
|C_{\psi(x)}\cap K_n|=|K_n|.
$$ 
\endproclaim

\flushpar
The map $f^{t_n-1}:K'_n\to K_n$ is monotone and differentiable hence we can
pull back the TI-pairs in $K_n$ into $K'_n$ and form the pairwise disjoint
collection: 
$$
\Cal{I}_n=\{I\subset K'_n|f^{t_n-1}(I) \text{ is a component of } 
                          C_{\psi(x)}\cap K_n\}.
$$  
From the previous lemma we get for every $I\in \Cal{I}_n$ an interval 
$T_I\subset K'_n$ such that $(T_I,I)$ is a TI-pair for $(V_x,U_x)$. Furthermore
$|\cup \Cal{I}_n|=|K'_n|$.

\flushpar
Let $\hat{\Cal{M}}_n$ consists of $I\in \Cal{I}_n$ with
\parindent=15pt
\item{1)} $I\cap [0,c_1]\ne \emptyset$;
\item{2)} $c_1\in T_I$.

\flushpar
From lemma 2.7 we know that the collection $\{T_I|I\in \hat{\Cal{M}}_n\}$
consists of nested intervals. In particular we can count 
$\hat{\Cal{M}}_n=\{\hat{M}_i|i\ge 0\}$  such that $T_{\hat{M}_{i+1}}\subset 
T_{\hat{M}_{i}}-\hat{M}_{i}$. Denote by $\hat{L}_i$ the component of 
$K'_n-\cup_{i\ge 1} \hat{M}_i$ next to $\hat{M}_i$ such that 
$c_1\notin [0,\hat{L}_i]$. 

\proclaim{Lemma 4.6} If $I\cap \hat{L}_i\ne \emptyset$ for some $i\ge 1$ and 
$I\in \Cal{I}_n$  then $T_I\subset \hat{L}_i$.

\flushpar
Furthermore there exists $\hat{\rho}>0$ such that 
$|\hat{L}_i|\ge \hat{\rho}|\hat{M}_i|$ for all $i\ge 0$.
\endproclaim

\demo{proof} The first statement follows directly from lemma 2.7.
The Koebe-Lemma gives $\hat{\rho}>0$ such that $|T_I|$ contains a 
$\hat{\rho}-$scaled neighborhood of $I$ where $(T_I,I)$ is a TI-pair for
$(V_ x,U_x)$. Again the statement follows from lemma 2.7: $\hat{L}_i$
contains a component of $T_{\hat{M}_i}-\hat{M}_i$.
\qed
\enddemo

\flushpar
Let $T_n$ be the symmetric interval $f^{-1}(K'_n)$. Observe that 
$\cup_{i\ge 0}\hat{L}_i\subset f(T_n)$. So define
$$
E_n=f^{-1}(\cup_{i\ge 0} \hat{L}_i).
$$

\flushpar
from the previous lemma we get

\proclaim{Lemma 4.7} Every $y\in E_n$ has a TI-pair for $(V_x,U_x)$. 
Furthermore this TI-pair is inside $T_n$.
\endproclaim

\flushpar
The next step is to show that $E_n\subset T_n$ has some universal metrical
properties. Denote the components of $T_n-E_n$ by 
$\Cal{M}_n=f^{-1}(\hat{\Cal{M}}_n)$. Again $\Cal{M}_n$ is countable, say 
$\Cal{M}_n=\{M_i|i\ge 0\}$. If $c\in \cup\Cal{M}_n$ then the corresponding 
component is $M_0$. If not  we define $M_0=\emptyset$. 

\flushpar
For $i\ge 1$ let $L_i$ be the component of $T_n-\cup_{i\ge 0} M_i$ next to 
$M_i$ but not between $M_i$ and $c$. Furthermore let $R_i$ be the component of
$T_n-M_i$ containing $c$. If $M_0\ne \emptyset$ define $R_0$ and $L_0$ to be
the components of $T_n-\cup_{i\ge 0} M_i$ next to $M_0$.

\proclaim{Lemma 4.8} There exists $\rho>0$ such that for $i\ge 0$
$$
|L_i|\ge \rho |M_i| \text{ and } |R_i|\ge \rho |M_i|.
$$
This bound is independent of $n\ge 1$.
\endproclaim

\demo{proof} 
For $i\ge 1$ $R_i$ contains one component of $T_n-\{c\}$. From the 
non-flatness of the critical and the fact that $T_n$ is a symmetric interval
we get $\rho_1>0$ with $|R_i|\ge \rho_1 |M_i|$.

\flushpar
Now $L_i\cup M_i$ is mapped monotonically onto $\hat{M}_j\cup \hat{L}_j$ for 
some $j\ge 1$. The space given by lemma 4.5 can be pulled back without being
distorted to much because the non-flatness of the critical point. Hence we get 
a constant $\rho_2>0$, depending on the behavior of the critical point such 
that
$$
|L_i|\ge \rho_2|M_i|.
$$ 

\flushpar
For $i=0$ the situation is even better: $L_0$ and $R_0$ are mapped onto some 
$\hat{L}_j$ and $M_0$ into the corresponding $\hat{M}_j$. Again the 
non-flatness  allow us to pull back space.
\qed
\enddemo

\demo{proof of proposition 4.3}
Let $\Cal R$ be the set of points whose orbit accumulates onto $c$. In section
2 we saw $|\Cal R|=1$. We are 
going to show that  the domain of the Markov map $M_x$ has a positive upper 
density in every point $y\in \Cal R$. From $|\Cal R|=1$ it follows that 
$|E_x|=1$.

\flushpar
Take $y\in \Cal R$ and $n\ge 1$. The first time the orbit of $y$ hits $T_n$ is
denoted by $s_n$. Let $H_n$ be the maximal interval containing $y$ such that
$f^{s_n}$ maps $H_n$ monotonically into $T_n$. Before we saw that there are no
renormalisations possible inside $U_x$. Hence the Lemma from the appendix
gives 
$$
f^{s_n}(H_n)=T_n.
$$  
Now let $A_n\subset H_n$ be the preimage of $E_n$. By using corollary 2.2 we 
get
a universal constant $\rho>0$ such that
$$
|A_n|\ge \rho |H_n|.
$$ 
The last observation to be made is $A_n\subset E_x$. The TI-pairs in $E_n$ are
part of $T_n$ hence they can be pulled back into $H_n$. So every point in 
$A_n$ has a TI-pair hence is in $E_x$.

\flushpar
Thus the upper density of $E_x$ in $y$ is bigger than $\rho$.  
\qed
\enddemo

\flushpar
Let us finish this section with the fundamental distortion property which holds
for Markov maps.

\flushpar
Fix $x\in \Cal N$ and consider the Markov map $M_x:E_x\to U_x$. Let the 
collection $\Cal{B}_0$ consists of the branches of $M_x$. That is it consists
of the components of $E_x$. Define inductively the collections of branches of 
$M_x^n$ as follows. An interval $I\subset U_x$ is in $\Cal{B}_{n+1}$ if 
$M_x(I)\in \Cal{B}_n$.

\proclaim{proposition 4.9} If $f$ satisfies the Markov-property, say
$M_x:E_x\to U_x$ has full domain, then there exists $K>0$ such that
$$
\frac{1}{K}\le \frac{|DM_x^n(y_1)|}{|DM_x^n(y_2)|}\le K
$$
for $y_1,y_2\in I\in \Cal{B}_n$ and  $n\ge 1$.
\endproclaim

\demo{Proof} Suppose that  $f$ also satisfies the Misiurewicz property, that
is, the critical point is not recurrent. In this case we may assume that the 
critical orbit does  not intersect $V_x$. Then it is easy to see that
the branches of $M_x^n$ have essentially bigger monotone extensions. Using
the Koebe-Lemma the proposition can be proved.

\flushpar
In the other case, the critical orbit is recurrent, we also find easily
monotone extensions. Indeed, using lemma 2.7 it is easy to see that every
branch of $M_x^n$ has a monotone extension to $V_x$. Again the Koebe-Lemma
finishes the proof.
\qed
\enddemo

\proclaim{Corollary 4.10} If $f$ satisfies the Markov-property then every 
0-dimensional closed invariant set has Lebesgue measure zero.
\endproclaim
 
\demo{Proof}
Let $X\subset [0,1]$ be a 0-dimensional closed invariant set and  let 
$y\in \Cal N$ be such that the Markov-Property holds for the transfer range
$(V_y,U_y)$. 

\flushpar
Suppose $|X|>0$. Then there exists a density point $x$ of $X$ with 
$x\in X\cap \Cal G$. Let $\{(T_n,I_n)\}_{n>0}$ be the TI-pairs for 
$(V_y,U_y)$ given by the Markov-Property. 
As in the proof of lemma 3.9 we show
$$
\lim_{n\to\infty} \frac{|X^c\cap U_x|}{|U_x|}
\le \lim_{n\to\infty} K \frac{|X^c\cap I_n|}{|I_n|}=0.
$$     
Hence, because $X$ is closed, $U_x\subset X$. Contradiction.
\qed
\enddemo

%% file: mark5.tex
\tolerance=3000

\bigskip
\centerline{\bf 5. The Ergodic Properties of Unimodal maps}
\bigskip

\flushpar
In this section the distortion results from the previous sections will be used
to give a measure theoretical description of $S-$unimodal dynamics. First
we will use the developed technics to prove two fundamental theorems which 
were already proved by Blokh and Lyubich in [BL1]: the ergodicity and the
Ergodic Classification Theorem. However we will include the relation with the
Markov-Property.  

\proclaim{Theorem 5.1} Every $S-$unimodal map which doesn't have periodic 
attractors is ergodic.
\endproclaim

\flushpar
The ergodicity follows directly from lemma 3.9.

\bigskip
\flushpar
We will fix in this section an $S-$unimodal map $f$ without periodic attractor.
Its critical point is denoted by $c$.

\bigskip
\flushpar 
Let $n\ge 0$ and $x\in [0,1]$. The maximal interval containing $x$ on which $f^n$
is monotone is denoted by $T_n(x)$. Furthermore the components of 
$T_n(x)-\{x\}$ are denoted by respectively $L_n(x)$ and $R_n(x)$. Define 
$r_n:[0,1]\to \bold{R}$ by
$$
r_n(x)=min\{|L_n(x)|,|R_n(x)|\}.
$$
This functions turn out to be fundamental for understanding the ergodic theory 
of unimodal maps.

\flushpar
Furthermore define
$$
r(x)=\limsup_{n\to\infty} r_n(x)  
$$
for every $x\in [0,1]$.

\flushpar
The ergodicity of $f$ implies the existence of a number $r\ge 0$ 
such that 
$$
r(x)=r   \text{ for a.e. } x\in [0,1].  
$$

\proclaim{Proposition 5.2} Let $f$ be an $S-$unimodal map without 
periodic attractor. The map $f$ has the Markov-Property if and only if $r>0$.
\endproclaim

\demo{Proof} The Markov-Property easily implies $r>0$.

\flushpar
Suppose we have $r>0$ for a map $f$ not satisfying the Markov-Property. Using 
the Contraction Principle we get that the function $\delta:[0,1]\to [0,1]$
$$
\delta(\epsilon)=\sup\{|I| | f^n \text{ maps the interval } I
                                 \text{ monotonically onto } T
                                 \text{ with }  |T|<\epsilon \}
$$ 
tends to zero as $\epsilon\to 0$.

\flushpar
Take $x\in\Cal{N}$ close enough to the critical point $c$ to get 
$\delta(\frac12 |V_x|)<\frac12 r$.

\flushpar 
Because $f$ doesn't satisfy the Markov-Property, we get $|E_x|<1$, which 
implies, by using the ergodicity, that for almost every point $y\in [0,1]$ 
the exists $n_y\ge 0$ such that for every $n\ge n_y$ with $f^n(y)\in U_x$  
$$
r_n(y)<\frac12 |V_x|. 
$$
Now take such a point $y$ and $n\ge n_y$. Let $s\ge 0$ be the smallest number 
such that $f^{n+s}(y)\in U_x$. Suppose that $r_{n+s}(y)$ is determined 
by the piece $R_{n+s}(y)$ of $T_{n+s}(y)$ (see definition $r_n(y)$).  
Now we claim that one piece, say $R_n(y)$ of $T_n(y)$ is mapped monotonically
onto $R_{n+s}(y)$. In fact this follows easily from the fact that $R_{n+s}(y)
\subset V_x$ and that $orb(\partial{U_x})\cap V_x=\emptyset$.

\flushpar
Then we conclude that $r_n(y)<\delta(\frac12 |V_x|)<\frac12 r$ for all 
$n\ge n_y$. Hence $\limsup r_n<\frac12 r$. Contradiction.
\qed
\enddemo

\bigskip
\flushpar
Using the number $r\ge 0$ we can give The Ergodic Classification of 
$S-$unimodal maps.
$$
\aligned
\Cal{P}&=\{f\in \Cal{U}| f \text{ has a periodic attractor }\}; \\
\Cal{C}&=\{f\in \Cal{U}| r=0\};\\
\Cal{I}&=\{f\in \Cal{U}| r>0\}. 
\endaligned
$$

\flushpar
Let us remember the definition given in [Mr] of an attractor (the ergodicity
of our maps makes the definition a bit simpler).

\proclaim{Definition 5.3} A closed invariant set $A\subset [0,1]$ is called an 
{\it attractor} if for almost every $x\in [0,1]$ 
$$
\omega(x)=A.
$$
\endproclaim

\proclaim{Theorem 5.4 (The Ergodic Classification Theorem)} Every $S-$unimodal 
map has an attractor $A$. It can be of three different types:
\parindent=15pt
\item{1)} if $f\in \Cal{P}$ then $A$ is a periodic orbit;
\item{2)} if $f\in \Cal{C}$ then $A=\omega(c)$ which is a minimal Cantor set;
\item{3)} if $f\in \Cal{I}$ then $A$ is the orbit of a periodic interval. 

\flushpar
In particular $f\in \Cal{I}$ if and only if $f$ has the Markov-Property.
\endproclaim

\demo{Proof} Consider an $S-$unimodal map $f$. If $f\in \Cal{P}$ we are 
finished.
\flushpar
If $f\in \Cal{C}$ then $f$ doesn't have the Markov-Property. 
Hence using Proposition 4.3 we get that
$\omega(c)$ is a minimal Cantor set. Now $r=0$ implies that almost all orbits
accumulates on $\omega(c)$. Again using the minimality of $\omega(c)$ the 
attractor turns out to be $\omega(c)$.

\flushpar
Take $f\in \Cal{I}$. By lemma 5.2
we get the Markov-Property for $f$, say $|E_x|=1$ for some $x\in \Cal{N}$. It 
easily follows from Proposition 4.9 that for almost every $x\in [0,1]$ 
$$
 U_x\subset \omega(x)
$$    
which implies $[c_2,c_1]\subset \omega(x)$. Because $\omega(x)\subset [c_2,c_1]$
for every $x\in [0,1]$ we get $A=[c_2,c_1]=\omega(x)$ for almost every 
$x\in [0,1]$.
\qed
\enddemo

\proclaim{Theorem 5.5} If $f$ is a $S-$unimodal map whose limit set 
$\omega(c)$ of the critical point is zero-dimensional then
$$
|\omega(c)|=0.
$$
In particular if $f\in \Cal{C}\cup \Cal{P}$ we have $|A|=0$.
\endproclaim 

\demo{Proof} Suppose that the critical point is not recurrent. Then there 
exists $x\in \Cal N$ with $\omega(c)\cap V_x=\emptyset$ or $\omega(c)
\subset \Lambda_x$. From lemma 2.5 we know $|\Lambda_x|=0$ which implies the
theorem.

\flushpar
The second case deals with the situation when $c\in \omega(c)$ but the limit 
set is not minimal.
In this case the map has the Markov-property. Corollary 4.10 tells that every
0-dimensional closed invariant set has Lebesgue measure zero. In particular
$|\omega(c)|=0$.
 
\flushpar
Now assume that $c\in \omega(c)$ and that $\omega(c)$ is minimal. Suppose $f$
is infinitely renormalizable. So $f|\omega(c)$ is injective. Together with 
lemma 3.9 we get $|\omega(c)|=0$. Hence we may assume that $f$ is not 
renormalizable. 

\flushpar 
Let $\{(V_{x_n},U_{x_n})\}_{n\ge 1}$ be the sequence of $\delta-$transfer 
ranges given by theorem 3.4. Fix $n\ge 1$. The minimality of $\omega(c)$ and 
lemma 2.5(3) implies that every point $x\in \omega(c)$ is contained in some 
component $I\subset C_{\psi(x_n)}$. So we get a finite collection of intervals 
$T_I$ covering $\omega(c)$. Here the intervals $T_I$ form together with the 
intervals $I$, $I\cap \omega(c)\ne \emptyset$, a TI-pair for 
$(V_{x_n},U_{x_n})$. Using the 
fact that these intervals $T_I$ are
nested, see lemma 2.7, and form a finite collection we can find a component
$I_n$ of $C_{\psi(x_n)}$ such that
\parindent=15pt
\item{1)} $I_n\cap \omega(c)\ne \emptyset$;
\item{2)} $(T_{I_n}-I_n)\cap \omega(c)=\emptyset$;

\flushpar
Because $(T_{I_n},I_n)$ is a TI-pair for the $\delta-$transfer range 
$(V_{x_n},U_{x_n})$ we get a $\rho_1>0$, not depending on $n\ge 1$ such that 
both components of $T_{I_n}-I_n$ are bigger than $\rho_1|I_n|$. Observe that
these components do not contain points of $\omega(c)$: we found some space in 
the Cantor set $\omega(c)$.

\flushpar
Let $t_n\ge 1$ be minimal such that $c_{t_n}\in I_n$. Using 2) above we get 
as before an interval $T'_n$ containing $c_1$ such that 
$f^{t_n-1}:T'_n\to T_{I_n}$ is an onto diffeomorphism. Let $T_n=f^{-1}(T'_n)$ 
and $M_n\subset T_n$ be the maximal interval with 
$f^{t_n}(M_n)\subset I_n$. Let $L_n$ and $R_n$ be the components of $T_n-M_n$. 
Using the fact that the critical point is non-flat we find $\rho>0$ such that
\parindent=15pt
\item{1)} $(T_n-M_n)\cap \omega(c)=\emptyset$;
\item{2)} $|L_n|=|R_n|\ge \rho|M_n|$.

\flushpar
Now we are able to prove that $\omega(c)$ does not have density points. For 
this let $x\in \omega(c)$. Because $c\in \omega(x)$ we can define $s_n\ge 0$ to
be the smallest number with $f^{s_n}(x)\in T_n$. 
Then apply the Lemma  from the appendix and we find for every $n\ge 1$ an interval $K_n$ around 
$x$ such that $f^{s_n}:K_n\to T_n$ is an onto diffeomorphism. Again the 
Contraction-principle assures that $|K_n|\to0$. Using the Macroscopic-Koebe-Lemma
and property 1) and 2) of $T_n$ above we see that $x$ is not a density 
point of $\omega(c)$. Hence $|\omega(c)|=0$. 
\qed
\enddemo

\proclaim{Theorem 5.7} Let $A$ be the topological attractor of $f\in \Cal I$.
Then $f|A$ is conservative.
\endproclaim

\demo{Proof} Let $D\subset A$ with positive Lebesgue measure. The map $f$ has 
the Markov-Property. Hence by Lemma 4.2 there exists $x\in \Cal N$ with
$|E_x|=1$. Because $D\subset A$ there exists $C\subset U_x$ with positive 
Lebesgue measure such that $f^n(C)\subset D$ for some $n\ge 0$. Let 
$D_0\subset D$ be the set of points whose orbits return to $D$.  Now by using
the Markov-Property we see that almost  every point of $D$ has a positive 
upper density for $D_0$. Hence $|D_0|=|D|$.
\qed
\enddemo

\flushpar
All the results are stated for $S-$unimodal maps. In fact it can be shown that 
the orbits of the intervals considered in the proofs of section 3
satisfy the necessary disjointness conditions needed for applying the 
$C^2-$Koebe-Lemma. Hence the results in this section also hold for 
$C^2-$unimodal maps. The difficulties for proving the $C^2$ versions of the 
Theorems in this 
paper occur when dealing with the orbits of TI-pairs. These orbits do not 
satisfy some disjointness conditions needed for applying the $C^2-$Koebe-Lemma.
However the intervals in these orbits are nested, as described by Lemma 2.7.
This property makes it possible to make the necessary estimates for 
applying the $C^2$-Koebe-Lemma. All results stated in this paper will turn out
to be true for $C^2$-unimodal maps.

%% file: marka.tex
\tolerance=3000

\bigskip
\centerline{\bf Appendix: Basic notions in real 1-dimensional dynamics}
\bigskip

\flushpar
Let $f:X\to X$ be a measurable map on the borel measure space $(X,\lambda)$.
We will give some basic definitions dealing with the ergodic theory of this
map.
\parindent=15pt
\item{0)} The orbit $\{x,f(x), f^2(x),...\}$ of a point $x$ is denoted by 
$orb(x)$ and the set of all limits of the $orb(x)$ is denoted by $\omega(c)$.
\item{1)} A borel set $A\subset X$ is invariant iff $f(X)\subset X$.
\item{2)} $f$ is called ergodic iff $X$ cannot be written as the union of
two disjoint invariant sets both with positive measure.
\item{3)} $f$ is called conservative iff for every set $D\subset X$ with 
$\lambda(D)>0$ the first return map $f_D$ on $D$ can be defined in almost 
every point of $D$.
\item{4)} A closed invariant set $D$ is called minimal if the orbit of every 
point in $D$ is dense in $D$.
\item{5)} acip stands for absolutely continuous invariant probability measure.
And acim stands for $\sigma-$finite  absolutely continuous invariant  measure.
 
\flushpar
Dealing with maps on the interval we will use the following notation:
\parindent=15pt
\item{6)} $\partial D$ is the boundary of the interval $D$;
\item{7)} The Lebesgue measure will be denoted by $|.|$;
\item{8)} A $\delta-$scaled neighborhood $T$ of the interval $I$ is an interval
such that both components of $T-I$ have length $\delta|I|$.  
\item{9)} Denote the maximal interval containing $x\in [0,1]$ on which $f^s$ is
monotone by $T_s(x)$. 

\bigskip
\flushpar
The collection $\Cal U$ consists of the $S-$unimodal maps. This are maps 
$f:[0,1]\to [0,1]$ having the following property:
\parindent=15pt
\item{1)} $f$ has negative Schwarzian derivative;
\item{2)} $f(0)=f(1)=0$;
\item{3)} there is exactly one point $c\in [0,1]$ where the derivative
of $f$ vanishes. Furthermore this {\it critical point} is non-flat: 
around $c$ the map behaves like $x\to x^{\alpha}$ ($\alpha>1$). The number
$\alpha>1$ is called the order of the critical point.

\bigskip
\flushpar
For every $S-$unimodal map $f$ the homeomorphism $\tau$ is defined to be the 
order reversing map satisfying $f\circ\tau=f$. The non-flatness of the critical 
point implies that this map is Lipschitz. Furthermore intervals of the form
$(x,\tau(x))$ are called symmetric.

\bigskip
\flushpar
An $S-$unimodal map $f$ is called renormalisable if there exists a symmetric 
interval $V$ such that the first return map to $V$ is of the form 
$f^n|V$, for some $n\ge 0$, and up to scaling $S-$unimodal. $f^n|V$ is 
called a renormalisation of $f$. It is called infinitely renormalisable if 
there are arbitrarily small symmetric intervals on which $f$ can be 
renormalised.

\flushpar
The limit set of the critical point of an infinitely renormalisable 
$S-$unimodal map is a minimal Cantor set. Furthermore the map acts like a 
homeomorphism on it.
 
\proclaim{Lemma} Let $f\in \Cal U$ be non-renormalisable and $V$ a symmetric
interval. If $s\ge 0$ is minimal such that $f^s(x)\in V$, $x\in [0,1]$, then
$V\subset f^s(T_s(x))$.
\endproclaim

\proclaim{Contraction Principle} Let $f\in \Cal U$ without periodic attractor.
Then for every $x\in [0,1]$
$$
|T_s(x)|\to 0
$$ 
when $s\to \infty$.
\endproclaim

\flushpar
The proofs of the above lemmas can be found for example in [MMS].

%% file: markr.tex
\tolerance=3000

\bigskip
\centerline{\bf References}
\bigskip

\parindent=20pt
\item{[BL1]} A.M.Blokh, M.Ju.Lyubich, {\it Attractors of maps of the interval},
Func. Ana. and Appl. vol. 21 (1987), 32-46.
\item{[BL2]} A.M.Blokh, M.Ju.Lyubich, {\it Measurable Dynamics of $S-$unimodal
maps of the interval}, preprint 1990/2 at SUNY.
\item{[G]} J. Guckenheimer, {\it Limit sets of $S-$unimodal maps with zero 
entropy}, Comm. Math. Phys. vol. 110, (1987), 133-160. 
\item{[GuJ]} J.Guckenheimer, S.Johnson, {\it Distortion of $S-$unimodal maps},
Ann. of Math. vol. 132 no. 1 (1990), 71-131.
\item{[M]} M.Martens, {\it Interval Dynamics}, Thesis at Technical University
of Delft, the Netherlands, (1990).
\item{[Mi]} M.Misiurewics, {\it Absolutely continuous measures for certain 
maps of the interval}, Publ. Math. IHES vol. 53 (1981) 17-51.
\item{[MMS]} M.Martens, W.de Melo, S.van Strien, {\it Julia-Fatou-Sullivan
theory for real one-dimensional dynamics}, to appear in Acta Math.  
\item{[Mr]} J.Milnor, {\it On the concept of attractors}, Commun. Math. Phys.
vol. 99 (1985) 177-195.